\documentclass[11pt]{article}
\usepackage{amssymb,psfig}
\usepackage{amsmath}

   \csname @addtoreset\endcsname{equation}{section}
    \csname @addtoreset\endcsname{figure}{section}
    \csname @addtoreset\endcsname{table}{section}

\newcounter{thanksnum}
\def\thanksnumber#1
{\setcounter{thanksnum}{\value{footnote}}\setcounter{footnote}{#1}%
                     \addtocounter{footnote}{-1}\footnotemark
                     \setcounter{footnote}{\value{thanksnum}}}
\def\newtheoremz#1{\@ifnextchar[{\@othmz{#1}}{\@nthmz{#1}}}

\def\@nthmz#1#2{%
\@ifnextchar[{\@xnthmz{#1}{#2}}{\@ynthmz{#1}{#2}}}

\def\@xnthmz#1#2[#3]{\expandafter\@ifdefinable\csname #1\endcsname
{\@definecounter{#1}\@addtoreset{#1}{#3}%
\expandafter\xdef\csname the#1\endcsname{\expandafter\noexpand
  \csname the#3\endcsname \@thmcountersepz \@thmcounterz{#1}}%
\global\@namedef{#1}{\@thmz{#1}{#2}}\global\@namedef{end#1}{\@endtheoremz}}}

\def\@ynthmz#1#2{\expandafter\@ifdefinable\csname #1\endcsname
{\@definecounter{#1}%
\expandafter\xdef\csname the#1\endcsname{\@thmcounterz{#1}}%
\global\@namedef{#1}{\@thm{#1}{#2}}\global\@namedef{end#1}{\@endtheoremz}}}

\def\@othmz#1[#2]#3{\expandafter\@ifdefinable\csname #1\endcsname
  {\global\@namedef{the#1}{\@nameuse{the#2}}%
\global\@namedef{#1}{\@thmz{#2}{#3}}%
\global\@namedef{end#1}{\@endtheoremz}}}

\def\@thmz#1#2{\refstepcounter
    {#1}\@ifnextchar[{\@ythmz{#1}{#2}}{\@xthmz{#1}{#2}}}

\def\@xthmz#1#2{\@begintheoremz{#2}{\csname the#1\endcsname}\ignorespaces}
\def\@ythmz#1#2[#3]{\@opargbegintheoremz{#2}{\csname
       the#1\endcsname}{#3}\ignorespaces}

\def\@thmcounterz#1{\noexpand\arabic{#1}}
\def\@thmcountersepz{.}
\def\@begintheoremz#1#2{ \trivlist \item[\hskip \labelsep{\bf #1\ #2}]}
\def\@opargbegintheoremz#1#2#3{ \trivlist
      \item[\hskip \labelsep{\bf #1\ #2\ (#3)}]}
\def\@endtheoremz{\endtrivlist}

\newtheorem{theorem}{Theorem}[section]

\newtheorem{proposition}{Proposition}[section]
\newtheorem{corollary}{Corollary}[section]

\newtheorem{remark}{Remark}[section]

\def\e{\varepsilon}

\def\o{\omega}
\def\O{\Omega}

\def\Y{{\cal Y}}
\def\F{{\cal F}}
\def\w{\widehat}
\def\Ind{{\,\rm Ind\,}}
\def\Ind{{\mathbb{I}}}

\def\R{{\bf R}}

\def\C{{\bf C}}

\def\X{{\cal X}}

\def\oo{\bar}

\def\U{{\cal U}}

\def\M{{\cal M}}

\newcommand{\be}{\begin{equation}}
\newcommand{\ee}{\end{equation}}
\newcommand{\bd}{\begin{displaymath}}
\newcommand{\ed}{\end{displaymath}}
\newcommand{\ba}{\begin{array}{ll}}
\newcommand{\ea}{\end{array}}
\newcommand{\baa}{\begin{eqnarray}}
\newcommand{\eaa}{\end{eqnarray}}
\newcommand{\baaa}{\begin{eqnarray*}}
\newcommand{\eaaa}{\end{eqnarray*}}
\font\sm=cmr10

\def\oo{\bar}

\def\Q{{\cal Q}}
\def\Re{{\rm Re\,}}

\def\U{{\cal U}}

\def\sinc{{\rm sinc\,}}
\date{Web-published:  29 November 2011. Revised: 10 August 2012 }
\title{
On  causal  band-limited  mean square
approximation}
\author{
Nikolai Dokuchaev\\
 {\sm Department of Mathematics \& Statistics, Curtin
University,}\\ {\sm  GPO Box U1987, Perth, 6845 Western
Australia}
}
\begin{document}
\maketitle \begin{abstract}  We study causal  dynamic approximation
of non-bandlimited processes by band-limited processes such that  a
part of the historical path of the underlying process  is
approximated in $L_2$-norm by the trace of a band-limited process.
This allows to cover the case of irregular non-smooth processes. We
show that this problem has an unique optimal solution. The
approximating band-limited process has unique extrapolation on
future times and can be interpreted as a optimal forecast.  To
accommodate the current flow of observations, the selection of this
band-limited process has to be changed dynamically. This can be
interpreted as a causal and linear filter that is not time
invariant.
\\
 {\bf Key words}: band-limited processes, causal filters,
sampling, low-pass filters,  prediction.
\\ AMS 2010 classification : 42A38, 
42B30, 
93E10 
\\
PACS 2008 numbers: 02.30.Mv, 
02.30.Nw,  
02.30.Yy, 
07.05.Mh,  
07.05.Kf 
\end{abstract}
\index{HIGHLIGHTS (FOR ONLINE VERSION) Causal  dynamic approximation
algorithm  of non-bandlimited processes by band-limited processes is
suggested. >  The approximation is sought for a  part of the
historical path of the underlying process by the trace of a
band-limited process >.  The approximating band-limited process is
obtained in time domain in a form of sinc series.}
\section{Introduction} We study causal dynamic approximation of
non-bandlimited processes by band-limited processes. It is known
that it is not possible to find an ideal low-pass causal linear
time-invariant filter. It is also known that the distance of the set
of these ideal low-pass time invariant filters from the set of all
causal filters is positive \cite{rema}. In addition, it is known
that optimal approximation of the ideal low-pass filter is not
feasible in the class of causal linear time-invariant filters (see,
e.g., \cite{D12} and references here). In the present paper, we are
trying to substitute the solution of these unsolvable problems by
solution of an easier problem where the filter is not necessary time
invariant. Our motivation is that, for some problems, time
invariancy for a filter is not crucial. For example, a typical
approach to forecasting in finance is to approximate the known path
of the stock price process by a smooth process that has an unique
extrapolation and accept this extrapolation as the forecast. This
procedure has to be done at current time; it is nor required that
the same forecasting rule will be applied at future times. We apply
this approach with the band-limited processes used as approximating
smooth predictable processes. More precisely, we suggest to
approximate in $L_2$-norm the known historical path of the process
by the trace of a band-limited process. In this setting, the
approximating curve does not necessary match the underlying process
at given sampling points. This is different from classical sampling
approach (see, e.g., \cite{jerry}). Similarly to
\cite{PFG}-\cite{PFG1}, our setting allows to cover the case of
irregular non-differentiable or discontinuous processes such as
historical stock prices in continuous time models. The difference is
that \cite{PFG}-\cite{PFG1} achieves point-wise matching for the
underlying process  being smoothed by a convolution operator; we
consider approximation of the underlying process directly using
different  methods. In \cite{PFG}-\cite{PFG1}, the estimate of the
error norm is given. In our setting, it is guaranteed that the
approximation generates  the error of the minimal norm.
\par
We show that an unique optimal solution of approximation problem
exits. The approximating process is derived in time domain in a form
of sinc series. To accommodate the current flow of observations, the
coefficients of these series and the related band-limited processes
have to be changed dynamically. It can be interpreted as a causal
and linear filter that is not time invariant.
\section{Definitions}
We denote by $L_2(D)$ the usual Hilbert space of complex valued
square integrable functions $x:D\to\C$, where $D$ is a domain.

For $x(\cdot)\in  L_2(\R)$, we denote by $X=\F x$ the function
defined on $i\R$ as the Fourier transform of $x(\cdot)$;
$$X(i\o)=(\F x)(i\o)= \int_{-\infty}^{\infty}e^{-i\o t}x(t)dt,\quad
\o\in\R.$$ Here $i=\sqrt{-1}$. For $x(\cdot)\in L_2(\R)$, the
Fourier transform $X$ is defined as an element of $L_2(\R)$ (more
precisely, $X(i\cdot)\in L_2(\R)$).

For a given $\O>0$, let $\U_{\O,\infty}=\{X(i\o)\in L_2(i\R):\
X(i\o)=0\ \hbox{for}\ |\o|>\O\}$, and  let  $\U_{\O,N}$ be  the set
of all $X\in U_{\O,\infty}$ such that there exists a sequence
$\{y_k\}_{k=-N}^N\in\C^{2N+1}$ such that
$X(i\o)=\sum_{k=-N}^{N}y_ke^{ik\o/\O}\Ind_{\{|\o|\le\O\}}$,  where
$\Ind$ is the indicator function.

For $N=+\infty$ and for integers $N\ge 0$, consider Hilbert spaces
$\Y_N$ such that $\Y_N=\C^{2N+1}$ for $N<+\infty$ and $\Y_N$ is the
set of all sequences $\{y_k\}_{k=-N}^N\in\C^{2N+1}$ such that
$\sum_{k=-\infty}^{\infty}|c_k|^2<+\infty$.

Let $s\in\R$ and $q<s$ be given; the case when  $q=-\infty$ is not
excluded. Consider Hilbert spaces of complex valued functions
$\X=L_2(-\infty,+\infty)$ and $\X_-=L_2(q,s)$.

Let $\O>0$ and $N$ be given (the case of $N=+\infty$ is not
excluded). Let $\X_{\O,N}$ be the subset of $\X_-$ consisting of
functions $x|_{(q,s]}$, where  $x\in\X$ are such that $x(t)=(\F^{-1}
X)(t)$ for $t\in[q,s]$ for some  $X(i\o)\in \U_{\O,N}$.
\begin{proposition}\label{propU}
For any $x\in\X_{\O,N}$, there exists an unique  $X\in\U_{\O ,N}$
such that $x(t)=(\F^{-1} X)(t)$ for $t\in[q,s]$.
\end{proposition}
\par
 For a Hilbert space $H$, we denote by
$(\cdot,\cdot)_{H}$ the corresponding inner product. We use notation
$\sinc(x)=\sin(x)/x$.
\section{Main results}
\subsection{Optimal band-limited approximation}
Let $x\in\X$ be a process. We assume that the path $x(s)|_{s\in
[q,s]}$ represents available historical data. Let Hermitian form
$F:\X_{\O,N}\times \X_-\to\R$ be defined as  \baaa F(\w
x,x)=\int_q^s|\w x(t)-x(t)|^2dt. \eaaa
\begin{theorem}\label{Th1} For
any $N\le +\infty$, there exists an unique solution  $\w x$  of the
minimization problem \baa \hbox{Minimize}\quad && F(\w
x,x)\quad\hbox{over}\quad \w x\in \X_{\O,N}.\label{min} \eaa
\end{theorem}
\begin{remark}\label{remF}
 By Proposition
\ref{propU}, there exists an unique extrapolation of the
band-limited solution $\w x(t)$  of problem (\ref{min}) on the
future time interval $(s,+\infty)$. It can be interpreted  as the
optimal forecast (optimal given $\O$ and $N$).
\end{remark}
\subsection{Optimal sinc coefficients}
To solve problem (\ref{min}) numerically, it is convenient to expand
$X(i\o)$ via Fourier series.

For a given $\O>0$, consider the mapping $\Q: \Y_N\to \X_{\O,N}$
such that $x=\Q y$ is such that $x(t)=(\F^{-1}X)(t)$ for a.e.
$t\in(q,s]$, where \baaa
X(i\o)=\sum_{t=-N}^{N}y_te^{it\o/\O}\Ind_{\{|\o|\le\O\}}.\eaaa
Clearly, this mapping is linear and continuous.
\par
Let Hermitian form $G:\Y_{N}\times \X_-\to\R$ be defined as \baa
G(y,x)=F(\Q y,x)=\int_q^s|\w x(t)-x(t)|^2dt, \quad \w x=\Q y.
\label{formG}\eaa
\begin{corollary}\label{corr1} For any $N\le +\infty$,
there exists an unique solution  $y$  of the minimization problem
\baa \hbox{Minimize}\quad && G(y,x)\quad\hbox{over}\quad y\in
\Y_{N}.\label{minY} \eaa
\end{corollary}
\par
Problem (\ref{min}) can be solved via problem (\ref{minY}); its
solution with $N<+\infty$ can be found numerically.

 \index{, and
$X_+=L_2(s,+\infty)ne nado?$ Let $\X_n$ be the subset of $\X_-$
consisting of functions $x\in\X$ such that $x(t)=(\F^{-1} X)(t)$ for
$t\in[q,s]$, where $X(i\o)\in \U_n$, and where $\U_n=\{X(i\o)\in
L_2(i\R):\ X(i\o)=0\ \hbox{for}\ |\o|\le \O\}$. Let
$x_n=x|_{[q,s]}-\w x$.
 Assume that $q=-\infty$ and that $x=\w
 x_c+\w x_n$, where $x_c\in \X_\O$ and $\w x_n\in \X_n$. }
\par
\subsection{Solution of problem (\ref{minY})} Let $N$ be given,  let
$Z$ be the set of all integers $z$ such that $|z|\le N$ if
$N<+\infty$, and let $Z$ be the set of all integers if $N=+\infty$.
Let \baaa X(i\o)=\sum_{k\in Z}y_ke^{ik\o\pi/\O}\Ind_{\{|\o|\le\O\}},
\eaaa where $\{y_k\}\in \Y_N$. Let $\w x=\F^{-1}X$. We have that
\baaa \w x(t)=\frac{1}{2\pi} \int_{-\O}^{\O}\left(\sum_{k\in Z}y_k
e^{ik\o\pi/\O}\right)e^{i\o t}d\o=\frac{1}{2\pi}
\sum_{k\in Z}y_k\int_{-\O}^{\O}e^{ik\o\pi/\O+i\o t}d\o\\
=\frac{1}{2\pi}\sum_{k\in Z}y_k \frac{e^{ik\pi+i\O t}- e^{-ik\pi-i\O
t}}{ik\pi/\O+it} =\frac{\O}{\pi}\sum_{k\in Z}y_k \sinc(k\pi+\O
t).\eaaa
\begin{remark}{\rm
Let $t[k]=-k\pi/\O$. Clearly, $\w x=\F^{-1}X$ is such that $\w
x(t[k])=y_k\cdot \O /\pi$, i.e., $y_k=\w x(t[k])\cdot \pi /\O$, and,
therefore,
 \baaa \w
x(t)=\sum_{k\in Z}\w x(t[k])\sinc(k\pi+\O t).\eaaa It gives
celebrated Sampling Theorem; see, e.g., \cite{jerry}. }\end{remark}
\begin{remark}{\rm We consider a setting when only the part
$x(t)|_{t\in[q,s]}$ of the path of the process is available at
current time $s<+\infty$. In this setting, sampling theorem is not
applicable. Our approximation can be considered as a modification of
the truncated sinc approximation (see, e.g., \cite{Jagerman},
\cite{jerry}). The difference is that the increasing of $N$ is not
related to extension the time interval $[q,s]$ in our setting. }
\end{remark}
\par
We have that \baa G(y,x)=\int_q^s|\w
x(t)-x(t)|^2dt=\int_q^s\left|\frac{\O}{\pi}\sum_{k\in Z}y_k
\sinc(k\pi+\O t)-x(t)\right|^2dt\nonumber\\=
(y,Ry)_{\Y_N}-2\Re(y,rx)_{\X_-}+(\rho x,x)_{\X_-}. \label{G}\eaa
Here $R:\Y_N\times \Y_N\to \Y_N$ is a linear bounded Hermitian
operator, $r:\X_-\to \Y_N$ is a bounded linear operator,
$\rho:\X_-\times \X_-\to \X_-$ is a linear bounded Hermitian
operator.
\par
It follows from the definitions that the operator $R$ is
non-negatively defined (it suffices to substitute $x(t)\equiv 0$
into the Hermitian form).

\subsection{The case when $N<+\infty$}
Up to the end of this paper, we assume that $N<+\infty$. In this
case, the space $\Y_N$  is finite dimensional, it follows that the
operator $R$ can be represented via a matrix $R=\{R_{km}\}\in
\C^{2N+1,2N+1}$, where $R_{km}=\oo R_{mk}$ and $(Ry)_k=\sum_{k=-N}^N
R_{km}y_m$.
\begin{theorem}
\label{ThP}
\begin{itemize}\item[(i)] For any $N<+\infty$, the operator $R$ is positively defined.
\item[(ii)] Problem (\ref{minY}) has a unique solution $\w y=R^{-1} rx$.
\item[(iii)] The components of the matrix $R$  can be found from the
equality \baa R_{km}=\index{
\frac{1}{(2\pi)^2}\int_q^s\overline{\left(\frac{e^{im\pi+i\O t}-
e^{-im\pi-i\O t}}{im\pi/\O+i t}\right)} \left( \frac{e^{ik\pi+i\O
t}- e^{-ik\pi-\O t}}{ik\pi/\O+i t}\right)dt\\=\frac{\O^2}{(2\pi)^2}
\int_q^s\frac{e^{-im\pi-i\O t}- e^{imt\O t}}{-im\pi-i\O t}\cdot
\frac{e^{ik\pi+i\O t}- e^{-ik\pi-\O t}}{ik\pi+i\O t}dt
\\= \frac{\O^2}{(2\pi)^2}\int_q^s\frac{-2i\sin(m\pi+\O t)}{-im\pi-i\O t}\cdot
\frac{2i\sin(k\pi+\O t)}{ik\pi+i\O t}dt\\=
\frac{\O^2}{\pi^2}\int_q^s\frac{\sin(m\pi+\O t)}{m\pi+\O t}\cdot
\frac{\sin(k\pi+\O t)}{k\pi+\O t}dt\\=}
\frac{\O^2}{\pi^2}\int_q^s\sinc(m\pi+\O t)\sinc(k\pi+\O
t)dt.\label{R}\eaa
\item[(iv)] The components of the vector  $rx=\{(rx)_k\}_{k=-N}^N$ can be found from the equality
\baa (rx)_{k}= \index{\int_q^s\overline{\left(\frac{e^{ik\pi+i\O t}-
e^{-ik\pi-i\O t}}{ik\pi/\O+i t}\right)} x(t) dt =\frac{\O}{\pi}
\int_q^s\frac{\sin(k\pi+\O t)}{k\pi+\O t}x(t) dt\\=}\frac{\O}{\pi}
\int_q^s\sinc(k\pi+\O t)x(t) dt. \label{r}\eaa
\end{itemize}
\end{theorem}
\begin{corollary} Let $\w y=\w  y(q,s)$ be the vector calculated as in Theorem \ref{ThP}, $\w y=\{\w
y_k\}_{k=-N}^N$. The process \baaa \w x(t)=\w x
(t,q,s)=\frac{\O}{\pi}\sum_{k\in Z}y_k \sinc(k\pi+\O t) \eaaa
\end{corollary}
represents the output of a causal filter that is linear but not time
invariant.
\section{Numerical experiments}
In the numerical experiments described below, we have used MATLAB
symbolic integration for calculation of integrals (\ref{R}) and
(\ref{r}) . The experiments show that some eigenvalues of $R$ are
quite close to zero. Because of the integration errors, some
eigenvalues of the calculated matrix $R$ are actually fluctuating
around zero despite the fact that, by Theorem \ref{ThP}, $R>0$.
Respectively, the error $E=\|R\w y-rx\|_{L_2(q,s)}|$ for the MATLAB
solution of the equation $R\w y=rx$ does  not vanish. This error
depends on the error tolerance parameter {\em tol} of MATLAB
integration operator $QUAD$ that was used; the default value is
$tol=10^{-6}$; we used $tol=10^{-8}$. Further, in our experiments,
we found that the error $E$ can be decreased by the replacing $R$ in
the equation $\w x=R^{-1}rx$ by $R_\e=R+\e I$, where $I$ is the unit
matrix and where $\e>0$ is small. In particular, for $\e=0.001$, the
corresponding error $E(\e)=\|R_\e^{-1}rx-\w y\|_{L_2(q,s)}<
\|R^{-1}rx-\w y\|_{L_2(q,s)}$, i.e., the approximation on $[q,s]$ is
better for $\w y=R_\e^{-1}rx$ calculated for $\e=0.001$ than for $\w
y=R^{-1}rx$ calculated for $\e=0$.

\begin{figure}[ht]
\centerline{\psfig{figure=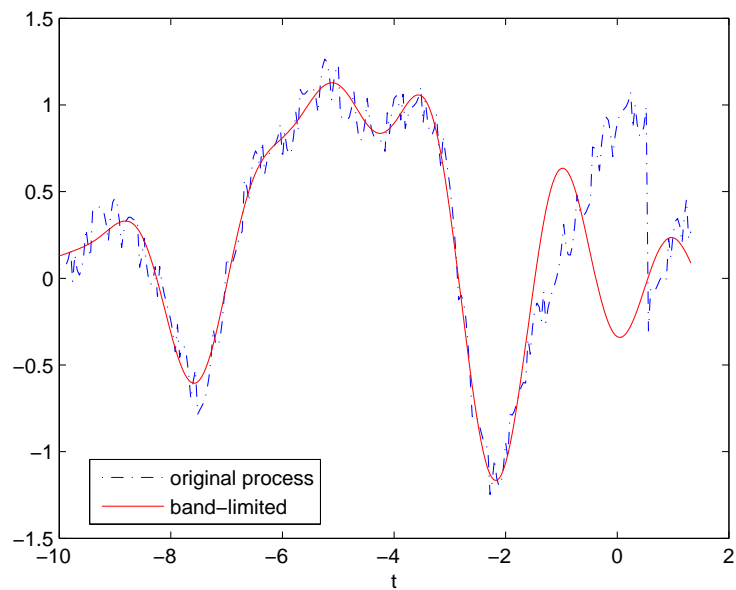,height=8.5cm}}
\caption[]{\sm Example of  $x(t)$ and band-limited process $\w x(t)$
approximating $x(t)$ on $(q,s]=(-12,-2]$, with $\O=4$, and $N=30$.}
\vspace{0cm}\label{fig-1}\end{figure}

\begin{figure}[ht]
\centerline{\psfig{figure=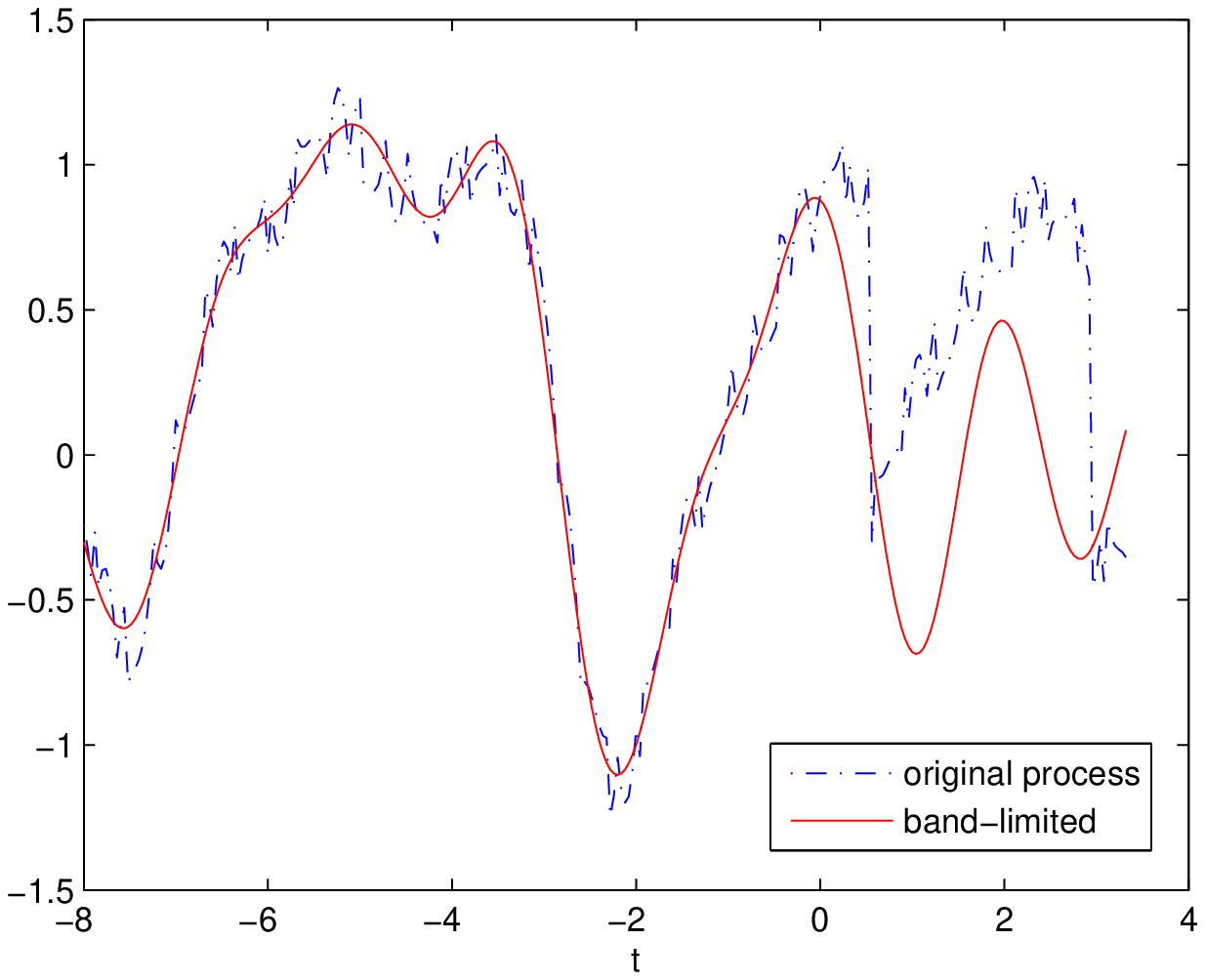,height=8.5cm}}
\caption[]{\sm Example of  $x(t)$ and band-limited process $\w x(t)$
approximating $x(t)$ on $(q,s]=(-10,0]$, with  $\O=4$, and $N=30$.}
\vspace{0cm}\label{fig0}\end{figure}

Figures \ref{fig-1} and \ref{fig0} show  examples of  a process
$x(t)$ and the band-limited process $\w x(t)$ approximating $x(t)$
on time intervals $(q,s]=(-12,-2]$ and $(q,s]=(-10,0]$,
respectively, calculated with $\e=0.001$ for $\O=4$ and $N=30$. As
expected, the change of the time interval from $(q,s]=(-12,-2]$ to
$(q,s]=(-10,0]$ results in the change of the approximating
band-limited process.

Note that the experiments demonstrate robustness with respect to the
changes of $N$. The curves of $\w x(t)$ will be almost the same if
we consider $N=50$ instead of $N=30$, when all other parameters are
the same. However, the error $E$ is larger for  large $N=100$, due
to accumulated larger error of integration.

\begin{figure}[ht]
\centerline{\psfig{figure=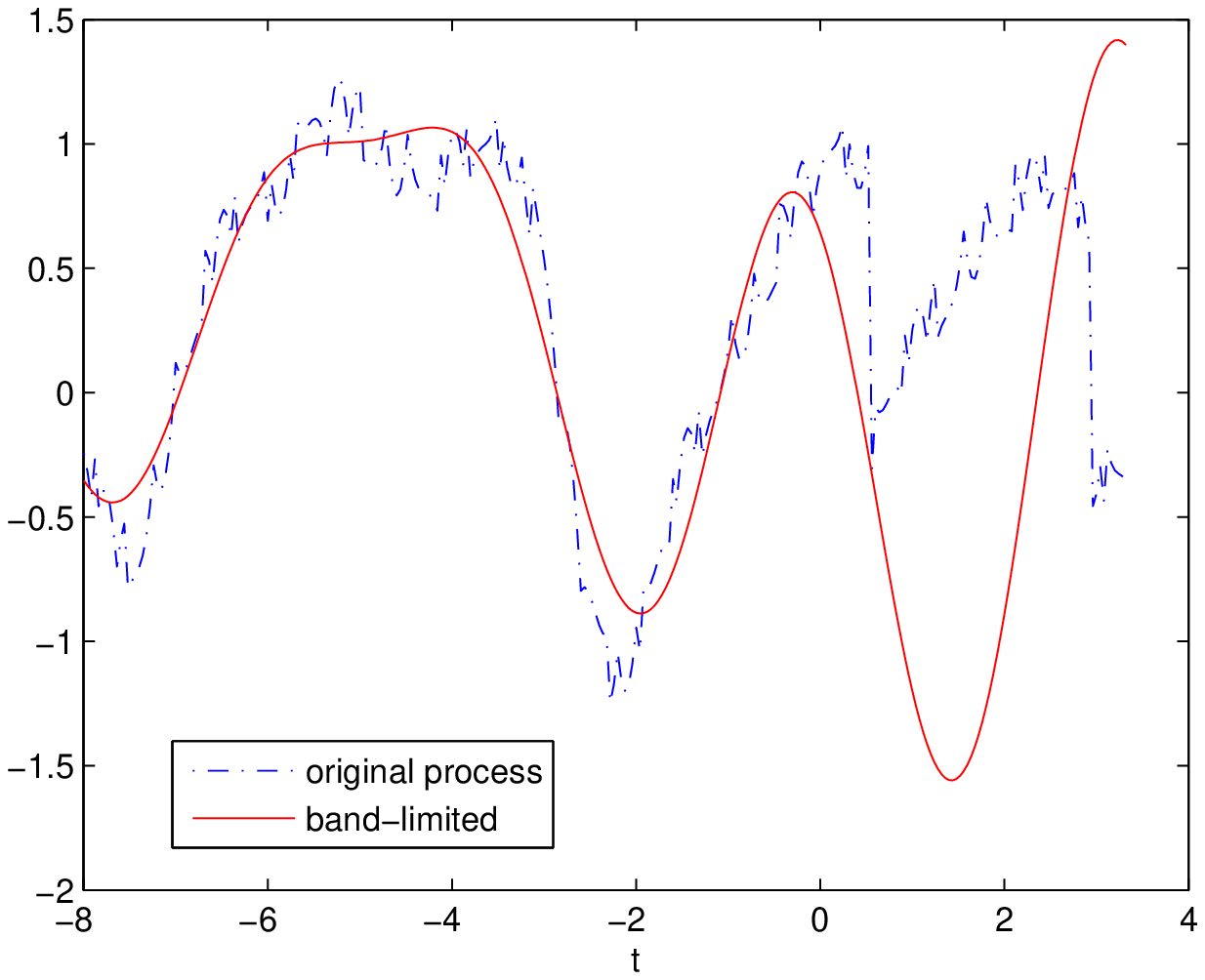,height=8.5cm}}
\caption[]{\sm Example of  $x(t)$ and band-limited process $\w x(t)$
approximating $x(t)$ on $(q,s]=(-10,0]$, with  $\O=2$, and $N=30$.}
\vspace{0cm}\label{fig2}\end{figure}

The shape of curves of $\w x(t)$ depends on the choice  $\O$. Figure
\ref{fig2} shows an example  of a process $x(t)$ and of the
band-limited process $\w x(t)$ approximating $x(t)$ on time interval
$(q,s]=(-10,0]$  calculated for $\O=2$, when all other parameters
are the same as for Figure \ref{fig0}.

By Remark \ref{remF}, the extrapolation of the process $\w x\in
\X_{\O,N}$ on the future time interval $(s,+\infty)$ can be
interpreted as the optimal forecast (optimal given $\O$ and $N$).
\begin{remark}\label{reme} {\rm
We have used the procedure of replacement $R$ by $R_\e=R+\e I$ with
small $\e>0$ to reduce the error  of calculation of the inverse
matrix for the matrix $R$ that is positively defined but is close to
a degenerate matrix. It can be noted that the same replacement could
lead to a meaningful setting  for the case when $\e>0$ is not small.
More precisely, it leads to optimization problem \baa
\hbox{Minimize}\quad &&
G(y,x)+\e^2\sum_{k=-N}^N|y_k|^2\quad\hbox{over}\quad y\in
\Y_{N}.\label{minYe} \eaa The solution restrains the norm of $y$,
and, respectively, the norm of $\w x$. }
\end{remark}
\par
\begin{figure}[ht]
\centerline{\psfig{figure=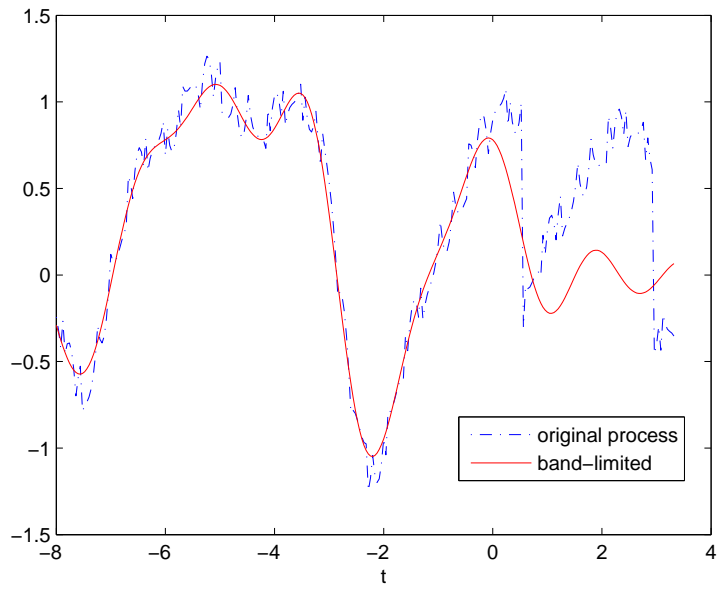,height=8.5cm}}
\caption[]{\sm Example of  $x(t)$ and band-limited process $\w x(t)$
calculated via solution of problem (\ref{minYe}) for $\e=0.05$,
$(q,s]=(-10,0]$, $\O=4$, and $N=30$.}
\vspace{0cm}\label{fig3}\end{figure}

Figure \ref{fig3} illustrates Remark \ref{reme} with  an example of
a process $x(t)$ and  the corresponding band-limited process $\w
x(t)$ calculated via solution of problem (\ref{minYe}) for
$\e=0.05$, when all other parameters are the same as for Figure
\ref{fig0}. This solution was obtained by replacement of $R$ by
$R_\e=R+\e I$ with $\e=0.05$.
\section{Appendix: proofs}
{\em Proof of Proposition \ref{propU}}.
  The  statement of this proposition  is known
in principle. It suffices to prove that if $x(\cdot)\in\X_{\O,N}$ is
such that $x(t)=0$ for $t\in(q,s]$, then $x(t)\equiv 0$. For the
sake of completeness,  we give below a proof.  For $C>0$, consider a
class $\M(C)$ of infinitely differentiable functions  $x(t):\R\to\R$
 such that there exists $M=M(x(\cdot))>0$ such that \baaa
&&\left\|\frac{d^kx}{dt^k}(\cdot)\right\|^2_{L_2(\R)}\le C^kM, \quad
k=0,1,2,.... \label{deriv1} \eaaa  Let $\M=\cup_{C>0}\M(C)$.  Any
$x\in \M$ is infinitely differentiable and
 such that there exists $C_1=C_1(x(\cdot))>0$ and $M_1=M_1(x(\cdot))>0$ such that \baaa
\sup_t\left|\frac{d^kx}{dt^k}(t)\right|\le C_1^kM_1. \eaaa Clearly,
$\X_{\O,N}\subset \M$. Therefore, any $x(\cdot)\in\X_{\O,N}$ is
analytic and allows the Taylor series expansion at any point with an
arbitrarily large radius of convergence. Consider the Taylor series
expansion at $t_0\in (q,s)$. Since all derivatives at this point are
equal to zero, the expansion is identically equal to zero.
 This
completes the proof of Proposition \ref{propU}. $\Box$
 \par
 {\em
Proof  of Theorem \ref{Th1}.} It suffices to  prove that $\X_{\O,N}$
is a closed linear subspace of $L_2(q,s)$.  In this case,  there
exists a unique projection $\w x$ of $x|_{[q,s]}$ on $\X_{\O,N}$,
and the theorem is proven.
\par
 Clearly, for  any $N\le +\infty$,
the set $U_{\O,N}$ is a closed linear subspace of $L_2(\R)$.
Consider a  mapping $Q:\U_{\O,N}\to \X_{\O,N}$ such that
$x(t)=(QX)(t)=(\F^{-1} X)(t)$ for $t\in[q,s]$. It is a linear
continuous operator. By Proposition \ref{propU}, it is a bijection.
 Since this mapping is continuous, it follows that
the inverse mapping $Q^{-1}: \X_{\O,N}\to U_{\O,N}$ is also
continuous (see Corollary in Ch.II.5 \cite{Yosida}, p.77). Since the
set $U_{\O,N}$ is a closed linear subspace of $L_2(\R)$, it follows
that $\X_{\O,N}$ is a closed linear subspace of $\X_-$. This
completes the proof of Theorem \ref{Th1}.
 $\Box$

 \par
{\em Proof  of Theorem \ref{ThP}.}    Let us prove statement (i). We
know that $R\ge 0$. Suppose that there exists $\oo y\in \C^{2N+1}$
such that $\oo y\neq 0$ and $R\oo y=0$. Let $r^*:\Y_N\to \X_-$ be
the adjoint operator to the operator $r^*:\X_-\to\Y_N$. If $r^*\oo
y\neq 0$ then there exists $x\in \X_-$ such that $G(\oo y,x)<0$,
which is not possible since $G(y,x)\ge 0$ for all $y,x$. Therefore,
$r^*\oo y=0$, i.e., $G(\oo y,x)=(\rho x,x)_{\X_-}$. Further, let $\w
y$ be a solution of problem (\ref{minY}). We have that  $G(\w
y,x)=G(\w y+\oo y,x)$. Hence $\w y+\oo y\neq \w y$ is another
solution of problem (\ref{minY}). This contradicts to Corollary
\ref{corr1} that states that this problem has an unique solution.
Statement (ii) follows from (i) and from classical theory of
quadratic forms. Statements (iii)-(iv) follow immediately from
representation (\ref{G}). This completes the proof of Theorem
\ref{ThP}. $\Box$
\subsection*{Acknowledgment}  This work  was supported by ARC grant of Australia DP120100928 to the author.

\end{document}